# EXPLICIT SOLUTIONS FOR THE SOLOMON-WILSON-ALEXIADES'S MUSHY ZONE MODEL WITH CONVECTIVE OR HEAT FLUX BOUNDARY CONDITIONS

**Domingo A. Tarzia** †‡[1]

*† Departamento de Matemática, FCE, Universidad Austral, Paraguay 1950,*

*S2000FZF  Rosario, Argentina.*

*‡ CONICET, Argentina.*

**Abstract**

We complete the Solomon-Wilson-Alexiades´s mushy zone model (Letters Heat Mass Transfer, 9 (1982), 319-324) for the one-phase Lamé-Clapeyron-Stefan problem. We obtain explicit solutions when a convective or heat flux boundary condition is imposed on the fixed face for a semi-infinite material. We also obtain the necessary and sufficient condition on data in order to get these explicit solutions. Moreover, when these conditions are satisfied the two problems are equivalents to the same problem with a temperature boundary condition on the fixed face and therefore an inequality for the coefficient which characterized one of the two free interfaces is also obtained.

**Nomenclature**

| | |
|---|---|
| $c$ | Specific heat, J/(kg ºC), |
| $-D_0 (<0)$ : | Temperature at the fixed face $x = 0$, ºC, |
| $-D_\infty (<0)$ | Bulk temperature at the fixed face $x = 0$, ºC, |
| $h_0 \ (h_0^*)$ | Coefficient that characterizes the transient heat transfer at x=0, kg/(Cºs^{5/2}), |
| $k$ | Thermal conductivity, W/(m ºC), |
| $\ell$ | Latent heat of fusion by unit of mass, J/kg, |
| $P_1$ | Phase-change process defined by conditions (1)-(6) |
| $P_2$ | Phase-change process defined by conditions (1)-(5) and (36) |
| $P_3$ | Phase-change process defined by conditions (1)-(5) and (43) |
| $q_0 \ (q_0^*)$ | Coefficient that characterizes the transient heat flux at x=0, kg/s^{5/2}, |
| $r = r(t) \ (> s(t))$ | Position of the liquid-mushy zone interface at time $t$, m, |
| $s = s(t)$ | Position of the solid-mushy zone interface at time $t$, m, |
| $t$ | Time, s, |
| $T$ | Temperature of the solid phase, ºC, |
| $x$ | Spatial coordinate, m, |

---

[1] Corresponding author. Tel.-fax: +54-341- 522-3093/3001
 E-mail address: DTarzia@austral.edu.ar

**Greek symbols**

| | |
|---|---|
| $\alpha = \dfrac{k}{\rho c}$ | Diffusivity coefficient, m$^2$/s, |
| $\gamma > 0$ | One of the two coefficients that characterizes the mushy zone, ºC, |
| $\varepsilon \in (0,1)$ | One of the two coefficients that characterizes the mushy zone, dimensionless, |
| $\nu\,(>\omega)$ | Coefficient that characterizes the free boundary $r(t)$ in Eq. (47), dimensionless, |
| $\rho$ | Density of mass, kg/m$^3$, |
| $\mu\,(>\xi)$ | Coefficient that characterizes the free boundary $r(t)$ in Eq. (10), dimensionless, |
| $\mu_\infty\,(>\xi_\infty)$ | Coefficient that characterizes the free boundary $r(t)$ in Eq. (28), dimensionless, |
| $\omega > 0$ | Coefficient that characterizes the free boundary $s(t)$ in Eq. (46), dimensionless |
| $\xi > 0$ | Coefficient that characterizes the free boundary $s(t)$ in Eq. (9), dimensionless |
| $\xi_\infty > 0$ | Coefficient that characterizes the free boundary $s(t)$ in Eq. (27), dimensionless |

**Key Words**

*Lamé-Clapeyron-Stefan Problem, PCM, Free boundary problem, Solomon-Wilson-Alexiades's mushy zone model, Explicit solutions, Convective boundary condition.*

## I. INTRODUCTION

Heat transfer problems with a phase-change such as melting and freezing have been studied in the last century due to their wide scientific and technological applications [1, 5-8, 11, 13, 18, 27]. A review of a long bibliography on moving and free boundary problems for phase-change materials (PCM) for the heat equation is shown in [24]. Explicit solutions to some free boundary problems was obtained in [3, 4, 9, 14 - 17, 19, 21, 28, 29]

We consider a semi-infinite material that is initially assumed to be liquid at its melting temperature which is assumed equals to 0 ºC. At time $t = 0$ a heat flux or a convective boundary condition is imposed at the fixed face $x = 0$, and a solidification process begins where three regions can be distinguished [20, 23]:

H1) liquid region at the temperature 0 ºC, in $x > r(t)$, $t > 0$;

H2) solid region at the temperature $T(x,t) < 0$, in $0 < x < s(t)$, $t > 0$ (with $s(t) < r(t)$);

H3) mushy region at the temperature $T(x,t) = 0$, in $s(t) \leq x \leq r(t)$, $t > 0$. The mushy region is considered isothermal and we make the following assumptions on its structure:

H3i) the material contains a fixed portion $\varepsilon\,\ell$ (with $0 < \varepsilon < 1$) of the total latent heat $\ell$ (see condition (3) below);

H3ii) the width of the mushy region is inversely proportional to the gradient of temperature (see condition (4) below).

Following the methodology given in [20, 23, 25] and the recent one in [26] we consider a convective boundary condition in Sections II to IV, and a heat flux condition in Sections V and VI at the fixed face $x = 0$ respectively. In both cases, we obtain explicit solutions for the temperature and the two free boundaries which define the mushy region. We also obtain, for both cases, the necessary and sufficient condition on data in order to get these explicit solutions given in Sections II and V respectively. Moreover, these two problems are equivalents to the same phase-change process with a temperature boundary condition on the fixed face $x = 0$ and therefore an inequality for the coefficient which characterized one of the two free interfaces is also obtained in Sections IV and VI. Moreover, in Section III we obtain the convergence of the phase-change process when the heat transfer coefficient goes to infinity.

## II. EXPLICIT SOLUTION WITH A CONVECTIVE BOUNDARY CONDITION

The phase-change process consists in finding the free boundaries $x = s(t)$ and $x = r(t)$, and the temperature $T = T(x,t)$ such that the following conditions must be verified (Problem ($P_1$)):

$$T_t - \alpha\, T_{xx} = 0\,, \quad 0 < x < s(t)\,, \quad t > 0 \, (\alpha = k\,/\,\rho c) \tag{1}$$

$$T\big(s(t),t\big) = 0\,, \quad t > 0 \tag{2}$$

$$kT_x(s(t),t) = \rho\ell[\varepsilon\dot{s}(t) + (1-\varepsilon)\dot{r}(t)]\,, \quad t > 0\,; \tag{3}$$

$$T_x(s(t),t)(r(t)-s(t)) = \gamma > 0\,, \quad t > 0 \,(\text{with } \gamma > 0)\,. \tag{4}$$

$$s(0) = r(0) = 0 \tag{5}$$

$$kT_x(0,t) = \frac{h_0}{\sqrt{t}}\big(T(0,t)+D_\infty\big)\,, \quad t > 0 \quad (h_0 > 0, D_\infty > 0)\,. \tag{6}$$

Condition (6) represents a convective boundary condition (Robin condition) at the fixed face $x = 0$ [2, 10, 12] with a heat transfer coefficient which is inversely proportional to the square root of the time [22, 25, 26, 30].

**Theorem 1.** *If the coefficient $h_0$ satisfies the inequality*

$$h_0 > \frac{1}{D_\infty}\sqrt{\frac{\gamma(1-\varepsilon)\rho\ell k}{2}} = h_0^* \tag{7}$$

*then the solution of problem* (1)-(6) *is given by:*

$$T(x,t) = -\frac{\frac{h_0 D_\infty \sqrt{\pi\alpha}}{k} erf(\xi)}{1 + \frac{h_0 \sqrt{\pi\alpha}}{k} erf(\xi)} \left[ 1 - \frac{erf\left(\frac{x}{2\sqrt{\alpha t}}\right)}{erf(\xi)} \right], \quad 0 < x < s(t), \ t > 0, \quad (8)$$

$$s(t) = 2\xi\sqrt{\alpha t}, \quad t > 0, \quad (9)$$

$$r(t) = 2\mu\sqrt{\alpha t}, \quad t > 0, \quad (10)$$

*with*

$$\mu = \xi + \frac{\gamma k}{2 D_\infty h_0 \sqrt{\alpha}} e^{\xi^2} \left[ 1 + \frac{h_0 \sqrt{\pi\alpha}}{k} erf(\xi) \right], \quad (11)$$

*and the coefficient* $\xi$ *is given as the unique solution of the equation:*

$$\frac{D_\infty c}{\ell \sqrt{\pi}} F(x) = G(x), \quad x > 0, \quad (12)$$

*where the real functions* $G$ *and* $F$ *are defined by:*

$$F(x) = \frac{e^{-x^2}}{\frac{k}{h_0 \sqrt{\pi\alpha}} + erf(x)} \quad, \quad G(x) = x + \frac{\gamma(1-\varepsilon)\sqrt{\pi}}{2 D_\infty} \frac{1}{F(x)} \quad, \quad x > 0. \quad (13)$$

**Proof.** Taking into account that $erf\left(\frac{x}{2\sqrt{\alpha t}}\right)$ is a solution of the heat equation (3) [6]

we propose as a solution of problem (1)-(6) the following expression:

$$T(x,t) = C_1 + C_2 \, erf\left(\frac{x}{2\sqrt{\alpha t}}\right), \quad 0 < x < s(t), \ t > 0, \quad (14)$$

where the two coefficients $C_1$ and $C_2$ must to be determined.

From condition (4) we deduce the expression (9) for the free boundary $s(t)$,

where the coefficient $\xi$ must be determined. From conditions (6) and (2) we deduce the

system of equations:

$$C_2 = \frac{h_0 \sqrt{\pi\alpha}}{k} \left( C_1 + D_\infty \right), \quad (15)$$

$$C_1 + C_2 \, erf(\xi) = 0, \quad (16)$$

whose solution is given by:

$$C_1 = -\frac{\frac{h_0 \sqrt{\pi\alpha}}{k} D_\infty erf(\xi)}{1 + \frac{h_0 \sqrt{\pi\alpha}}{k} erf(\xi)}, \quad C_2 = \frac{h_0 D_\infty \sqrt{\pi\alpha}}{k} \frac{1}{1 + \frac{h_0 \sqrt{\pi\alpha}}{k} erf(\xi)}, \quad (17)$$

and then we get expression (8) for the temperature.

From condition (4) we deduce expression (10) for the interface $r(t)$ and expression (11) for $\mu$. From condition (3) we deduce equation (12) for the coefficient $\xi$. Functions $F_3$ and $G$ have the following properties:

$$F(0^+) = \frac{h_0 \sqrt{\pi \alpha}}{k} > 0, \quad F(+\infty) = 0^+, \quad F'(x) < 0, \quad \forall x > 0, \tag{18}$$

$$G(0^+) = \frac{(1-\varepsilon)\gamma k}{2 D_\infty h_0 \sqrt{\alpha}} > 0, \quad G(+\infty) = +\infty, \quad G'(x) > 0, \quad \forall x > 0. \tag{19}$$

Therefore, we deduce that equation (12) has a unique solution when the coefficient $h_0$ satisfies the inequality

$$\frac{D_\infty c}{\ell \sqrt{\pi}} F(0^+) > G(0^+) \quad \Leftrightarrow \quad h_0^2 > \frac{\gamma(1-\varepsilon)\rho \ell k}{2 D_\infty^2}, \tag{20}$$

i.e. inequality (7) holds. □

**Corollary 2.** *If the coefficient* $h_0$ *satisfies inequality* (7) *then the temperature, defined by* (8), *verifies the following inequalities:*

$$-D_\infty < T(0,t) \leq T(x,t) < 0, \quad 0 < x < s(t), \ t > 0. \tag{21}$$

**Proof.** From (8) we obtain:

$$T(0,t) = -\frac{\dfrac{h_0 D_\infty \sqrt{\pi \alpha}}{k}}{1 + \dfrac{h_0 \sqrt{\pi \alpha}}{k} erf(\xi)} = -\frac{D_\infty}{1 + \dfrac{k}{h_0 \sqrt{\pi \alpha} \ erf(\xi)}} > -D_\infty, \quad \forall t > 0. \tag{22}$$

Moreover, from (8) and (22) we also get

$$\begin{aligned}
T(x,t) + D_\infty &= \frac{D_\infty}{1 + \dfrac{h_0 \sqrt{\pi \alpha}}{k} erf(\xi)} \left[ 1 + \frac{h_0 \sqrt{\pi \alpha}}{k} erf\left( \frac{x}{2\sqrt{\alpha t}} \right) \right] \\
&\geq \frac{D_\infty}{1 + \dfrac{h_0 \sqrt{\pi \alpha}}{k} erf(\xi)} = T(0,t) + D_\infty > 0, \quad 0 < x < s(t), \ t > 0
\end{aligned} \tag{23}$$

that is (21) holds. □

## III. ASYMPTOTIC BEHAVIOR WHEN THE COEFICIENT $h_0 \to +\infty$

Now, we will obtain the asymptotic behaviour of the solution (8)-(12) of problem (1)-(6) when the heat transfer coefficient is large, that is when $h_0 \to +\infty$.

For any coefficient $h_0$ satisfying inequality (7) we will denote the temperature $T$ and the two free boundaries $s$ and $r$ by $T = T(x, t, h_0)$, $x = s(t, h_0)$ and $x = r(t, h_0)$ respectively, with coefficients $\xi = \xi(h_0)$ and $\mu = \mu(h_0)$. We will also denote with $F(x, h_0)$ and $G(x, h_0)$ the functions defined in (13). We have the following result:

**Theorem 3.** *We obtain the following limits:*

$$\lim_{h_0 \to \infty} T(x, t, h_0) = T_\infty(x, t), \quad \lim_{h_0 \to \infty} s(t, h_0) = s_\infty(t), \quad \lim_{h_0 \to \infty} r(t, h_0) = r_\infty(t), \quad (24)$$

*where $T_\infty(x,t), s_\infty(t)$ and $r_\infty(t)$ are the solutions of the following phase-change process with mushy region:* (1)-(5) *and*

$$T(0, t) = -D_\infty, \quad t > 0, \quad (25)$$

*instead of the boundary condition* (6).

**Proof.** The solution of problem (1)-(5) and (25) is given by [20]:

$$T_\infty(x, t) = -D_\infty \left[ 1 - \frac{erf\left(\dfrac{x}{2\sqrt{\alpha t}}\right)}{erf(\xi_\infty)} \right], \quad 0 < x < s_\infty(t), \ t > 0, \quad (26)$$

$$s_\infty(t) = 2\xi_\infty \sqrt{\alpha t}, \quad t > 0, \quad (27)$$

$$r_\infty(t) = 2\mu_\infty \sqrt{\alpha t}, \quad t > 0, \quad (28)$$

with

$$\mu_\infty = \xi_\infty + \frac{\gamma \sqrt{\pi}}{2 D_\infty} e^{\xi_\infty^2} erf(\xi_\infty), \quad (29)$$

and the coefficient $\xi_\infty$ given as the unique solution of the equation:

$$G_1(x) = \frac{D_\infty c}{\ell \sqrt{\pi}}, \quad x > 0, \quad (30)$$

where the real function $G_1$ is defined by:

$$G_1(x) = \frac{G_\infty(x)}{F_\infty(x)}, \quad x > 0. \quad (31)$$

with

$$G_\infty(x) = \left[ x + \frac{\gamma(1 - \varepsilon)\sqrt{\pi}}{2 D_\infty} \frac{1}{F_\infty(x)} \right] = \lim_{h_0 \to \infty} G(x, h_0), \quad x > 0. \quad (32)$$

$$F_\infty(x) = \frac{e^{-x^2}}{erf(x)} = \lim_{h_0 \to \infty} F(x, h_0), \quad x > 0 \,. \tag{33}$$

Then,

$$\lim_{h_0 \to \infty} \xi(h_0) = \xi_\infty, \quad \lim_{h_0 \to \infty} \mu(h_0) = \mu_\infty \,, \tag{34}$$

and therefore, the limits (24) hold.

**Remark 1.** By studying the real functions $F(x, h_0)$ and $G(x, h_0)$ we can obtain the order of the convergence:

$$0 < \xi_\infty - \xi(h_0) = O\left(\frac{1}{h_0}\right) \text{ when } h_0 \to \infty \,. \tag{35}$$

## IV  EQUIVALENCE BETWEEN THE MUSHY ZONE MODELS WITH CONVECTIVE AND TEMPERATURE BOUNDARY CONDITIONS

We consider the problem ($P_2$) defined by the conditions (1) – (5) and temperature boundary condition

$$T(0, t) = -D_0 < 0 \,, \quad t > 0 \,, \tag{36}$$

at the fixed face $x = 0$, whose solution was given in [20]. We have the following property:

**Theorem 4.** *If the coefficient $h_0$ satisfies inequality (7) then Problem ($P_1$), defined by conditions (1)-(6), is equivalent to Problem ($P_2$), defined by conditions (1)-(5) and (36), when the parameter $D_0$ in Problem ($P_2$) is related to parameters $h_0$ and $D_0$ in Problem ($P_1$) by the following expression:*

$$D_0 = \frac{D_\infty \, erf(\xi)}{\dfrac{k}{h_0 \sqrt{\pi \alpha}} + erf(\xi)} > 0 \tag{37}$$

*where the coefficient $\xi$ is given as the unique solution of equation (12) for Problem ($P_1$) or as the unique solution of equation:*

$$G_2(x) = \frac{D_0 \, c}{\ell \sqrt{\pi}} \,, \quad x > 0 \,, \tag{38}$$

*for Problem ($P_2$) where the real function $G_2$ is defined by:*

$$G_2(x) = \frac{G_0(x)}{F_\infty(x)}, \quad G_0(x) = \left[ x + \frac{\gamma(1-\varepsilon)\sqrt{\pi}}{2D_0} \frac{1}{F_\infty(x)} \right], \quad x > 0. \quad (39)$$

**Proof.** If the coefficient $h_0$ satisfies inequality (7) then the solution of the Problem ($P_1$) is given by $(8) - (12)$. Taking into account that:

$$T(0,t) = -\frac{\frac{h_0\sqrt{\pi\alpha}}{k}D_\infty erf(\xi)}{1 + \frac{h_0\sqrt{\pi\alpha}}{k}erf(\xi)} = -\frac{D_\infty erf(\xi)}{\frac{k}{h_0\sqrt{\pi\alpha}} + erf(\xi)} < 0, \quad t > 0 \quad (40)$$

then we can define the Problem ($P_2$) by imposing the temperature boundary condition (36) with data $D_0$ given in (37). By using this data $D_0$ in the Problem ($P_2$) and the method developed in [26] we can prove that the solutions of both Problems ($P_1$) and ($P_2$) are the same and then the two problems are equivalents. □

**Corollary 5.** *If the coefficient $h_0$ satisfies inequality* (7) *then the coefficient $\xi$ of the solid-mushy zone interface of Problem ($P_2$) verifies the following inequality:*

$$erf(\xi) < \frac{D_\infty D_0}{D_\infty - D_0}\sqrt{\frac{2c}{\pi\gamma(1-\varepsilon)\ell}}, \quad \forall D_\infty > D_0 . \quad (41)$$

*Then,*

$$erf(\xi) < D_0\sqrt{\frac{2c}{\pi\gamma(1-\varepsilon)\ell}} . \quad (42)$$

**Remark 2.** The real functions $G_\infty$, defined in (32), and $G_0$, defined in (39), are similar; the difference between them are the parameters $D_\infty$ or $D_0$ used in each definition.

## V. EXPLICIT SOLUTION WITH A HEAT FLUX BOUNDARY CONDITION

Now, we will consider a phase-change process which consists in finding the free boundaries $x = s(t)$ and $x = r(t)$, and the temperature $T = T(x,t)$ such that the following conditions must be verified (Problem ($P_3$)): conditions (1) - (5), and

$$kT_x(0,t) = \frac{q_0}{\sqrt{t}}, \quad t > 0 \quad (q_0 > 0). \quad (43)$$

Condition (43) represents the heat flux at the fixed face $x = 0$ characterized by a coefficient which is inversely proportional to the square root of the time [22].

**Theorem 6.** *If the coefficient $q_0$ satisfies the inequality*

$$q_0 > \sqrt{\frac{\gamma(1-\varepsilon)\rho\ell k}{2}} = q_0^* \tag{44}$$

*then the solution of problem* (1)-(5) *and* (43) *is given by:*

$$T(x,t) = -\frac{q_0\sqrt{\pi\alpha}\; erf(\omega)}{k}\left[1 - \frac{erf\left(\dfrac{x}{2\sqrt{\alpha t}}\right)}{erf(\omega)}\right] < 0, \quad 0 < x < s(t),\; t > 0, \tag{45}$$

$$s(t) = 2\omega\sqrt{\alpha t}, \quad t > 0, \tag{46}$$

$$r(t) = 2\nu\sqrt{\alpha t}, \quad t > 0, \tag{47}$$

*with*

$$\nu = \omega + \frac{\gamma k}{2q_0\sqrt{\alpha}}e^{\omega^2}, \tag{48}$$

*and the coefficient $\omega > 0$ given as the unique solution of the equation:*

$$G_3(x) = \frac{q_0}{\rho\ell\sqrt{\alpha}}, \quad x > 0, \tag{49}$$

*where the real function $G_3$ is defined by:*

$$G_3(x) = \left[x + \frac{\gamma(1-\varepsilon)k}{2q_0\sqrt{\alpha}}e^{x^2}\right]e^{x^2}, \quad x > 0. \tag{50}$$

**Proof.** Following the proof of the Theorem 1, we propose as a solution of problem (1)-(5) and (43) the following expression:

$$T(x,t) = A_1 + A_2\; erf\left(\frac{x}{2\sqrt{\alpha t}}\right), \quad 0 < x < s(t),\; t > 0, \tag{51}$$

where the two coefficients $A_1$ and $A_2$ must to be determined.

From condition (2) we deduce expression (46) for the free boundary $s(t)$, with the coefficient $\omega$ to be determined. From conditions (2) and (43) we deduce:

$$A_1 = -\frac{q_0\sqrt{\pi\alpha}}{k}erf(\omega), \quad A_2 = \frac{q_0\sqrt{\pi\alpha}}{k}, \tag{52}$$

and then we get expression (45) for the temperature.

From condition (4) we deduce expression (47) for the interface $r(t)$ and expression (48) for $v$. From condition (3) we deduce equation (49) for the coefficient $\omega$. Since function $G_3$ has the following properties:

$$G_3(0^+) = \frac{\gamma(1-\varepsilon)k}{2q_0\sqrt{\alpha}} > 0, \quad G_3(+\infty) = +\infty, \quad G_3^{'}(x) > 0, \quad \forall x > 0, \qquad (53)$$

we can deduce that equation (49) has a unique solution when the coefficient $q_0$ satisfies the inequality

$$\frac{q_0}{\rho\ell\sqrt{\alpha}} > G_3(0^+) \quad \Leftrightarrow \quad q_0^2 > \frac{\gamma(1-\varepsilon)\rho\ell k}{2}, \qquad (54)$$

i.e. inequality (44). □

**Remark 3.** We have a relationship between $q_0^*$ (the lower limit for the coefficient $q_0$ in order to have a phase-change process with a mushy region with a heat flux boundary condition at $x = 0$) and $h_0^*$ (the lower limit for the coefficient $h_0$ in order to have a phase-change process with a mushy region with a convective boundary condition at $x = 0$) given by:

$$q_0^* = D_\infty h_0^*. \qquad (55)$$

## VI EQUIVALENCE BETWEEN THE MUSHY ZONE MODELS WITH HEAT FLUX AND TEMPERATURE BOUNDARY CONDITIONS

Following Section IV, we will now study the relationship between the Problems ($P_3$) and ($P_2$). We have the following property:

**Theorem 7.** *If the coefficient $q_0$ satisfies inequality* (44) *then Problem* ($P_3$)*, defined by conditions* (1)-(5) *and* (43)*, is equivalent to Problem* ($P_2$)*, defined by conditions* (1)-(5) *and* (36)*, when the parameter $D_0$ in Problem* ($P_2$) *is related to the parameter $q_0$ in Problem* ($P_3$) *by the following expression:*

$$D_0 = \frac{q_0\sqrt{\pi\alpha}}{k} erf(\omega) > 0 \qquad (56)$$

*where the coefficient $\omega$ is given as the unique solution of equation* (49) *for Problem* ($P_3$) *or as the unique solution of equation* (38) *for Problem* ($P_2$).

**Proof.** If the coefficient $q_0$ satisfies inequality (44) then the solution of Problem ($P_3$) is given by (45) − (49). Taking into account that:

$$T(0,t) = -\frac{q_0\sqrt{\pi\alpha}}{k}\, erf(\omega) < 0, \quad t > 0, \tag{57}$$

we can define the Problem ($P_2$) by imposing the temperature boundary condition (36) with the data $D_0$ given in (56). By using this data $D_0$ in Problem ($P_2$) and the method developed in [26] we can prove that the solutions of both Problems ($P_3$) and ($P_2$) are the same and then the two problems are equivalents. □

**Corollary 8.** *If the coefficient $q_0$ satisfies inequality* (44) *then the coefficient $\xi$ of the solid-mushy zone interface of the Problem* ($P_2$) *verifies inequality* (42) *which is the same that we have obtained through the equivalence between Problems* ($P_1$) *and* ($P_2$).

By using the results of this work, we can now obtain new explicit expression for the determination of one or two unknown thermal coeffcient through a phase-change process with a mushy zone by imposing an overspecified convective boundary condition at the fixed face $x = 0$. This will complete and improve the results obtained previously in [23].

## CONCLUSIONS

The goal of this paper is to complete the Solomon-Wilson-Alexiades's model for a mushy zone model for phase-change materials when a convective boundary or a heat flux condition at the fixed face $x = 0$ is imposed. In both cases, explicit solutions for the temperature and the two free boundaries which define the mushy region was obtained; and, for both cases, the necessary and sufficient conditions on data in order to get these explicit solutions are also obtained. Moreover, the equivalence of these two phase-change process with the one with a temperature boundary condition on the fixed face $x = 0$ was obtained. On the other hand, the convergence of the phase-change process with mushy zone when the heat transfer coefficient goes to infinity was also obtained.

## AKNOWLEDGEMENTS

The present work has been partially sponsored by the Projects PIP No 0534 from CONICET - Univ. Austral, Rosario, Argentina, and AFOSR-SOARD Grant FA9550-14-1-0122.

# REFERENCES

[1] V. Alexiades, A.D. Solomon, Mathematical modeling of melting and freezing processes, Hemisphere-Taylor & Francis, Washington, 1996.

[2] P.M. Beckett, A note on surface heat transfer coefficients, Int. J. Heat Mass Transfer 34 (1991) 2165-2166.

[3] A.C. Briozzo, D.A. Tarzia, Explicit solution of a free-boundary problem for a nonlinear absorption model of mixed saturated-unsaturated flow, Adv. Water Resources, 21 (1998) 713-721.

[4] P. Boadbridge, Solution of a nonlinear absorption model of mixed saturated-unsaturated flow, Water Resources Research, 26 (1990) 2435-2443.

[5] J.R. Cannon, The one-dimensional heat equation, Addison-Wesley, Menlo Park, California, 1984.

[6] H.S. Carslaw, C.J. Jaeger, Conduction of heat in solids, Clarendon Press, Oxford, 1959.

[7] J. Crank, Free and moving boundary problem, Clarendon Press, Oxford, 1984.

[8] A. Fasano, Mathematical models of some diffusive processes with free boundary, MAT – Serie A, 11 (2005) 1-128.

[9] F. Font, S.L. Mitchell, T.G. Myers, One-dimensional solidification of supercooled melts, Int.J. Heat Mass Transfer, 62 (2013) 411-421.

[10] S.D. Foss, An approximate solution to the moving boundary problem associated with the freezing and melting of lake ice, A.I.Ch.E. Symposium Series, 74 (1978) 250-255.

[11] S.C. Gupta, The classical Stefan problem. Basic concepts, modelling and analysis, Elsevier, Amsterdam, 2003.

[12] C.L. Huang, Y.P. Shih, Perturbation solution for planar solidification of a saturated liquid with convection at the hall, Int. J. Heat Mass Transfer, 18 (1975) 1481-1483.

[13] V.J. Lunardini, Heat transfer with freezing and thawing, Elsevier, London, 1991.

[14] T.G. Myers, F. Font, On the one-phase reduction of the Stefan problem with a variable phase change temperature, Int. Comm. Heat Mass Transfer, 61 (2015) 37-41.

[15] T.G. Myers, S.L. Mitchell, F. Font, Energy conservation in the one-phase supercooled Stefan problem, Int. Comm. Heat Mass Transfer, 39 (2012) 1522-1525.

[16] M.F. Natale, D.A. Tarzia, "Explicit solutions to the two-phase Stefan problem for Storm's type materials", J. Physics A: Mathematical and General, 33 (2000) 395-404.

[17] C. Rogers, On a class of reciprocal Stefan moving boundary problems, Z. Angew. Math. Phys., (2015), DOI 10.1007/s00033-015-0506-1.

[18] L.I. Rubinstein, The Stefan problem, American Mathematical Society, Providence, 1971.

[19] N.N. Salva, D.A. Tarzia, "Explicit solution for a Stefan problem with variable latent heat and constant heat flux boundary conditions", J. Math. Analysis Appl., 379 (2011) 240-244.

[20] A.D. Solomon, D.G. Wilson, V.A. Alexiades, A mushy zone model with an exact solution, Letters Heat Mass Transfer, 9 (1982) 319-324.

[21] A.D. Solomon, D.G. Wilson, V. Alexiades, Explicit solutions to change problems, Quart. Appl. Math., 41 (1983) 237-243.

[22] D.A. Tarzia, An inequality for the coefficient $\sigma$ of the free boundary $s(t) = 2\sigma\sqrt{t}$ of the Neumann solution for the two-phase Stefan problem, Quart. Appl. Math., 39 (1981) 491-497.

[23] D.A. Tarzia, Determination of unknown thermal coefficients of a semi-infinite material for the one-phase Lamé-Clapeyron (Stefan) problem through the Solomon-Wilson-Alexiades mushy zone model, Int. Comm. Heat Mass Transfer, 14 (1987) 219-228.

[24] D.A. Tarzia, A bibliography on moving-free boundary problems for heat diffusion equation. The Stefan problem, MAT - Serie A 2 (2000) 1-297.

[25] D.A. Tarzia, An explicit solution for a two-phase unidimensional Stefan problem with a convective boundary condition at the fixed face, MAT – Serie A, 8 (2004) 21-27.

[26] D.A. Tarzia, Relationship between Neumann solutions for two-phase Lamé-Clapeyron-Stefan problems with convective and temperature boundary conditions, Thermal Science, (2015), In Press. See arXiv-1406.0552.

[27] A.B. Tayler, Mathematical models in applied mechanics, Clarendon Press, Oxford, 1986.

[28] V.R. Voller, F. Falcini, Two exact solutions of a Stefan problem with varying diffusivity, Int. J. Heat Mass Transfer, 58 (2013) 80-85.

[29] V. R. Voller, J.B. Swenson, C. Paola, An analytical solution for a Stefan problem with variable latent heat, Int. J. Heat Mass Transfer, 47 (2004) 5387-5390.

[30] S.M. Zubair, M.A. Chaudhry, Exact solutions of solid-liquid phase-change heat transfer when subjected to convective boundary conditions, Heat Mass Transfer, 30 (1994) 77-81.